\documentclass[12pt]{article}

\usepackage{latexsym}
\usepackage{amsmath}
\usepackage{amsthm}
\usepackage{amsfonts}
\usepackage{amsgen}
\usepackage{pictex}
\usepackage{amssymb}
\usepackage{wasysym}
\usepackage{authblk}
\usepackage{titlefoot}
\usepackage[T1]{fontenc}
\usepackage{eulervm} 
\usepackage{tikz}
\usetikzlibrary{spy}
\usetikzlibrary{matrix,arrows}
\usetikzlibrary{topaths}
  \usepackage{graphicx}
\usetikzlibrary{fadings}
\usetikzlibrary{patterns}

\usepackage{graphicx}
\usepackage{wrapfig}

\usepackage{enumitem}
\usepackage[perpage,symbol*]{footmisc}
\usepackage[font=footnotesize,labelfont=bf]{caption}
\usepackage[utf8]{inputenc}
\input amssym.def 
\input amssym.tex 
\newcommand{\Proof}{\noindent{\sc Proof} \hspace{0.2in}}

\newcommand{\smallbox}[1]{\leavevmode\thinspace\hbox{\vrule\vtop{\vbox 
   {\hrule\kern1pt\hbox{\strut\thinspace{#1}\thinspace}} 
   \kern1pt\hrule}\vrule}\thinspace} 

\define\peg{{}^\frown\!}

\define\qed{\hfill $\square$}

\newcount\skewfactor
\def\mathunderaccent#1#2 {\let\theaccent#1\skewfactor#2
\mathpalette\putaccentunder}
\def\putaccentunder#1#2{\oalign{$#1#2$\crcr\hidewidth
\vbox to.2ex{\hbox{$#1\skew\skewfactor\theaccent{}$}\vss}\hidewidth}}

\newcount\skewfactor
\def\mathunderaccent#1#2 {\let\theaccent#1\skewfactor#2
\mathpalette\putaccentunder}
\def\putaccentunder#1#2{\oalign{$#1#2$\crcr\hidewidth
\vbox to.2ex{\hbox{$#1\skew\skewfactor\theaccent{}$}\vss}\hidewidth}}

\newcommand{\cK}{{\cal K}}

\newcommand{\lesdot}{\mathrel{\mathord{<}\!\!\raise 0.8
pt\hbox{$\scriptstyle\circ$}}}

   \def\nofork^#1_#2{\putforkinmargin
        \unionstick^{\textstyle #1}_{\textstyle #2}}

\newbox\noforkbox \newdimen\forklinewidth
\setbox0\hbox{$\textstyle\bigcup$}
\setbox1\hbox to \wd0{\hfil\vrule width \forklinewidth depth \dp0
                        height \ht0 \hfil}
\setbox\noforkbox\hbox{\box1\box0\relax}
\def\unionstick{\mathop{\copy\noforkbox}\limits}
\def\nonfork#1#2_#3{#1\unionstick_{\textstyle #3}#2}
\def\nonforkin#1#2_#3^#4{#1\unionstick_{\textstyle #3}^{\textstyle #4}#2}

\def\doesforkin_#1^#2{\biguplus\limits_{\textstyle #1}^{\textstyle
#2}}
\newtheorem{theorem}{Theorem}[section] 
\newtheorem*{theorem*}{Theorem}

\newtheorem{claim}[theorem]{Claim}

\newtheorem{notation}[theorem]{Notation}

\title{Infinitary Logics and Abstract Elementary Classes}

\author[1]{Saharon Shelah}
\author[2]{Andrés Villaveces}

\affil[1]{Hebrew University of Jerusalem / Rutgers University}
\affil[2]{Universidad Nacional de Colombia - Bogotá}
\begin{document} 
 
\maketitle

\unmarkedfntext{Research of both authors partially supported by NSF grant no: DMS 1833363.}
\unmarkedfntext{Research of first author partially supported by Israel Science Foundation (ISF)
grant no: 1838/19.}

\abstract{
    We prove that every abstract elementary class (a.e.c.) with LST number $\kappa$ and vocabulary
$\tau$ of cardinality $\leq \kappa$ can be axiomatized in the logic ${\mathbb
L}_{\beth_2(\kappa)^{+++},\kappa^+}(\tau)$. In this logic an a.e.c. is therefore an EC class
rather than merely a PC class. This constitutes a major improvement on the level of definability
previously given by the Presentation Theorem. 
As part of our proof, we define the \emph{canonical tree} $\mathcal S={\mathcal
S_\cK}$ of an a.e.c. $\cK$. This turns out to be an interesting combinatorial object of the class, beyond the aim of
our theorem.
Furthermore, we study a connection between the
sentences defining an a.e.c. and the relatively new infinitary logic $L^1_\lambda$.}


\section*{Introduction}

Given an abstract elementary class (a.e.c.) $\cK$, in vocabulary $\tau$ of size $\leq
\kappa=LST(\cK)$, we do two main things:
\begin{itemize}
    \item We provide an infinitary sentence \emph{in the same vocabulary} $\tau$ of the a.e.c. that
	axiomatizes $\cK$.
    \item We also provide a version of the ``Tarski-Vaught-criterion,'' adapted to a.e.c.'s: when
	$M_1\subseteq M_2$, for $M_1,M_2\in \cK$, we will provide necessary and sufficient
	\emph{syntactic} conditions for $M_1\prec_\cK M_2$. These will depend on a certain sentence
	holding \emph{only} in $M_2$.
\end{itemize}

The two proofs hinge on a new combinatorial object: a \textbf{canonical
tree} ${\mathcal S}_\cK$ for an a.e.c. $\cK$.
${\mathcal S}_\cK$ is a well-founded tree of models in $\cK$, all of them of
cardinality equal to $LST(\cK)$. The tree ${\mathcal S}_\cK$ encodes all possible embeddings between
models of size $LST(\cK)$ in $\cK$. Although it is a purely combinatorial object, it encapsulates
enough information on the a.e.c. $\cK$ as to enable us to axiomatize it.

The Presentation Theorem~\cite[p. 424]{Sh:88} is central to the development of stability for abstract
elementary classes: notably, it enables Ehrenfeucht-Mostowski techniques for classes that have
large enough models. This has as an almost immediate consequence stability below a
categoricity cardinal and opens the possibility of a relatively advanced
classification/stability theory in that wider setting.

The Presentation Theorem had provided a way to capture an a.e.c. as a PC-class: by expanding
its vocabulary with infinitely many function symbols, an a.e.c. may be
axiomatized by an infinitary formula. Although for the stability-theoretical
applications mentioned this expansion is quite useful, the question as to whether it is
possible to axiomatize an a.e.c. with an infinitary sentence \emph{in the same
vocabulary} of the a.e.c. is natural. Here we provide a positive
solution: given an a.e.c $\cK$ 
we provide an infinitary sentence \emph{in the same original vocabulary}
$\varphi_\cK$  whose models are exactly those in $\cK$. Therefore, unlike the
situation in the Presentation Theorem, here the class turns out to be an EC Class, not a PC
class.

The main idea is that a ``canonical tree of models'', each of size the LST-number of the class,
the tree of height $\omega$ ends up providing enough tools; the sentence essentially describes
all possible maps from elements of this tree into arbitrary potential models in the class. A
combinatorial device (a partition theorem theorem on well-founded trees due to
Komjáth and Shelah~\cite{KoSh:796}) is
necessary for our proof.


The two main theorems:

\begin{theorem*}[Theorem~\ref{aecDefinableInfinitary}]
    \textbf{(Axiomatization of an a.e.c. in $\tau$ by an infinitary sentence in $\tau$.)}
    Let $\kappa=LST(\cK)+|\tau|$ for an abstract elementary class $\cK$ in vocabulary $\tau$,
    and let $\lambda=\beth_2(\kappa)^{+ +}$.
    Then
    there is a sentence $\psi_\cK$ in the logic ${\mathbb L}_{\lambda^+,\kappa^+}(\tau)$ such that
    $\cK=Mod(\psi_\cK)$.
\end{theorem*}

Our second theorem provides a syntactic characterization of being a $\prec_\cK$-elementary submodel,
a kind of syntactic``Tarski-Vaught'' criterion for AECs. The precise statement is lengthy and
requires notation built later here. But here is a version.

\begin{theorem*}[Description of the main point of Theorem~\ref{newTVlike}]
    \textbf{(A syntactic ``Tarski-Vaught'' criterion for $\prec_\cK$-elementarity.)}
    If $M_1\subseteq M_2$ are $\tau=\tau_\cK$-structures, then the following are equivalent:
    \begin{itemize}
	\item $M_1\prec_\cK M_2$,
	\item given \textbf{any} tuple $\bar{a}\in M_1$ of length $\kappa$, we may find a
	    $\tau$-structure $N_{\bar{a}}$ of size $\kappa$ and an isomorphism $f_{\bar{a}}$ from
	    $N_{\bar{a}}$ onto some $N_{\bar{a}}^*\subseteq_\tau M_2$ such that $\bar{a}\subseteq
	    N_{\bar{a}}^*$ and $M_2$ satisfies a formula (called $\varphi_{N_{\bar{a}},\lambda+1,1}$;
	    it will be part of a hierarchy of formulas we'll define), on the elements of
	    $N_{\bar{a}}^*$. The formula $\varphi_{\dots}$ is in the same logic where we can
	    axiomatize $\cK$; namely, ${\mathbb L}_{\lambda^+,\kappa^+}(\tau)$.
    \end{itemize}
\end{theorem*}

The second part of the previous theorem, the characterization of being $\prec_\cK$-elementary, thus
amounts to the following: for every tuple in $M_1$, the model $M_2$ satisfies a formula describing
the fact that the tuple may be covered by another tuple that has the ``eventual tree extendibility''
property described by the formula $\varphi_{\dots}$ that we define in the next section.

We finish the paper by connecting our axiomatization with logics close to Shelah's logic
$L^1_\kappa$ and other logics similar to it, recently studied by D\v{z}amonja and
Väänänen~\cite{DzVa} and by Väänänen and Villaveces~\cite{VaVi}.

We wish to thank Xavier Caicedo, Mirna D\v{z}amonja, Rami Grossberg, Juliette Kennedy and Jouko Väänänen, as well as
the referee of this paper, for useful
comments and remarks on earlier versions of this paper.

\section{Canonical trees and sentences for a.e.c.'s}

We fix an a.e.c. $\cK$ for the remainder of this paper, with vocabulary $\tau$. We also fix
$\kappa=LST(\cK)\geq |\tau|$ and we let $\lambda=\beth_2(\kappa)^{++}$.
Without loss of generality we assume that all models in $\cK$ are of cardinality $\geq
\kappa$. Furthermore, we will use for the sake of convenience an ``empty model'' called $M_{\rm empt}$ with the property that
$M_{\rm empt}\prec_\cK M$ for all $M\in \cK$.

\subsection{The canonical tree of an a.e.c.}

We now build a canonical object for our abstract elementary class $\cK$, ${\mathcal S}={\mathcal
S}_\cK$. This will be a tree with $\omega$-many levels, consisting of models in $\cK$ of size
$\kappa$, organized in a way we now describe. To prove our results, we will use the tree ${\mathcal
S}_\cK$ to ``test'' membership in $\cK$ and ``depths'' of possible extensions.

\begin{notation}
We fix the following notation for the rest of this paper.
\begin{itemize}
    \item We first fix a sequence of (different) elements $(a^*_\alpha\mid
	    \alpha<\kappa\cdot
	\omega)$ in some model in $\cK$.
    \item $\bar{x}_n:=\langle x_\alpha\mid \alpha<\kappa\cdot n\rangle$,
    \item $\bar{x}_{=n}:= \langle x_\alpha\mid \alpha\in \left\{ \kappa\cdot n + \zeta \mid
	    \zeta < \kappa\right\}\rangle$.
	\item $\bar{z}_{[\kappa]}$ denotes $\langle z_\alpha\mid \alpha<\kappa\rangle$.

\end{itemize}   
\end{notation}

\noindent
We now define the \textbf{canonical tree}\label{canonicaltree} of $\cK$:

\begin{itemize}
    \item ${\mathcal S}_n:=\big\{ M\in \cK\mid M$ has universe $\left( a^*_\alpha
	\right)_{\alpha <\kappa\cdot n}$  and $m<n$ implies $M\restriction\left( a^*_\alpha
    \right)_{\alpha <\kappa\cdot m}\prec_\cK M \big\}$,
    \item ${\mathcal S}={\mathcal S}_\cK:=\bigcup_{n<\omega}{\mathcal S}_n$; this is a tree with $\omega$
	levels under $\prec_\cK$ (equivalenty under $\subseteq$, by our definition of each level).
\end{itemize}

We use this tree in our proof to test properties of the class $\cK$. The key point about ${\mathcal
S}_\cK$ is that it contains information not just on models in the class of cardinality
$\kappa=LST(\cK)$ but more importantly on the way they embed into one another.

%

\subsection{Formulas and sentences attached to $\cK$}

We now define by induction on $\gamma<\lambda^+$ formulas
\[ \varphi_{M,\gamma,n}(\bar{x}_n),\]
for every $n$ and $M\in {\mathcal S}_n$ (when $n=0$ we may omit $M$). \textbf{We build all these
formulas within the logic ${\mathbb L}_{\lambda^+,\kappa^+}(\tau)$.}

\begin{description}
    \item[Case 1]: $\gamma=0$\\ If $n=0$ then the formula $\varphi_{0,0}$ is
	$\top$ (the sentence denoting ``truth''). Assume $n>0$. Then
	\[ \varphi_{M,0,n}(\bar{x}_n):=\bigwedge {\rm Diag}^n_\kappa(M),\]
	where ${\rm Diag}^n_\kappa(M)$ is the set $\Big\{ \varphi(x_{\alpha_0},\dots,x_{\alpha_{k-1}})\mid
	\alpha_0,\dots,\alpha_{k-1}<\kappa\cdot n$,
    $\varphi(y_0,\dots,y_{k-1})$ is an atomic or a negation of an atomic formula and $M\models 
	\varphi(a^*_{\alpha_0},\dots,a^*_{\alpha_{k-1}})\Big\}$.
    \item[Case 2]: $\gamma$ a limit ordinal\\ Then
	\[
	\varphi_{M,\gamma,n}(\bar{x}_n):=\bigwedge_{\beta<\gamma}\varphi_{M,\beta,n}(\bar{x}_n).\]
    \item[Case 3]: $\gamma=\beta+1$\\
	Let $\varphi_{M,\gamma,n}(\bar{x}_n)$ be the formula
	\[\forall
	    \bar{z}_{[\kappa]}\bigvee_{\stackrel{N\succ_\cK M}{N\in {\mathcal S}_{n+1}}} \exists
	    \bar{x}_{=n}\left[
	    \varphi_{N,\beta,n+1}(\bar{x}_{n+1})\wedge 
    \bigwedge_{\alpha<\kappa}\bigvee_{\delta< \kappa\cdot(n+1)}z_\alpha=x_\delta \right]\]
\end{description}

By construction, all these formulas $\varphi_{M,\gamma,n}(\overline{x}_n)$ belong to ${\mathbb
L}_{\lambda^+,\kappa^+}(\tau)$. When $n=0$, for all $\gamma<\lambda^+$, these
formulas are really \emph{sentences} $\varphi_{M_{\rm empty},\gamma,0}$ (as $\overline{x}_0$ is
the empty sequence). As said above, in this case we just denote them by $\varphi_{\gamma,0}$. These sentences may be
understood as ``approximations'' to a sentence fully characterizing the a.e.c. $\cK$.

Our first aim is to prove that the sentences $\varphi_{\gamma,0}$ indeed end up
axiomatizing the a.e.c. $\mathcal{K}$, when $\gamma$ is large enough.

%

\section{Characterizing $\cK$ by its canonical sentence}

In this section we prove the first main theorem:

\begin{theorem}
    There is a sentence $\psi_\cK$ in the logic ${\mathbb L}_{\lambda^+,\kappa^+}(\tau)$ such that
    $\cK=Mod(\psi_\cK)$.
    \label{aecDefinableInfinitary}
\end{theorem}

Our first aim in this section is to prove that every model $M\in
\cK$ satisfies $\varphi_{\gamma,0}$, for all $\gamma<\lambda^+$.

In order to achieve this, we prove the following (more elaborate) statement, by induction on
$\gamma$. (Recall $\varphi_{\gamma,0}$ is the sentence $\varphi_{M_{\rm empty},\gamma,0}$ for
each $\gamma<\lambda^+$.)

\begin{claim}\label{inductSat}
    Given $\gamma<\lambda^+$, $M\in \cK$, $n<\omega$, $N\in {\mathcal S}_n$, $f:N\to M$ a
    $\prec_\cK$-embedding (if
    $n=0$, $f$ is empty) \underline{then} $M\models \varphi_{N,\gamma,n}[\langle
    f(a^*_\alpha\mid \alpha<\kappa\cdot n)\rangle]$.
\end{claim}

Before starting the proof, notice that in the statement of the Claim, when $n=0$, we have that
$f$ is empty and $\varphi_{\gamma,0}$ is a sentence. Notice also as $\gamma$ grows, the sentences
$\varphi_{\gamma,0}$ capture ever more involved properties of the model $M$. Thus, when $\gamma=0$,
$\varphi_{0,0}$ holds trivially; for $\gamma=1$, $M\models \varphi_{1,0}$ means $M$ satisfies
$\forall \bar{z}_{[\kappa]}\bigvee_{\stackrel{N\succ_\cK M_{empt}}{N\in {\mathcal S}_1}} \exists
\bar{x}_{=1}\left[ \varphi_{N,0,1}(\bar{x}_{1})\wedge \bigwedge_{\alpha<\kappa}\bigvee_{\delta<
\kappa\cdot 1}z_\alpha=x_\delta \right]$. This means that given any subset $Z\subseteq M$ of size at
    most $\kappa$, there is some $N\in {\mathcal S}_1$, the first level of the canonical tree, such
    that the image of $N$ under some embedding $f:N\to M$, $f(X)$, covers $X$. In short, this
    amounts to saying that $M$ is densely covered by images of models in $\cK$ of size $\kappa$.

    When $\gamma=2$, we know a bit more: parsing the sentence, $M\models \varphi_{2,0}$ means that
    in $M$,
    \[\forall \bar{z}_{[\kappa]}\bigvee_{\stackrel{N\succ_\cK M_{empt}}{N\in {\mathcal S}_1}} \exists
    \bar{x}_{=1}\left[ \varphi_{N,1,1}(\bar{x}_{1})\wedge \bigwedge_{\alpha<\kappa}\bigvee_{\delta<
\kappa\cdot 1}z_\alpha=x_\delta\right] .\]
Parsing again, this means that 
\[\forall
	\bar{z}_{[\kappa]}\bigvee_{\stackrel{N\succ_\cK M_{empt}}{N\in {\mathcal S}_1}} \exists
	\bar{x}_{=1}\left[ \forall \bar{z'}_{[\kappa]} \bigvee_{\stackrel{N'\succ_\cK N}{N'\in
	    {\mathcal S}_2}}\exists \bar{x}_{=2} \varphi_{N',0,2}(\bar{x}_{2})\wedge
    \bigwedge_{\alpha<\kappa}\bigvee_{\delta< \kappa\cdot 2}z'_\alpha=x_\delta\right.\]
    \[ \left. \qquad \qquad \wedge \bigwedge_{\alpha<\kappa}\bigvee_{\delta<\kappa\cdot
2}z_\alpha=x_\delta\right].\]

    What this long formula says is that given any subset $Z\subseteq M$ there is some $N$ in level
    $1$ of the tree ${\mathcal S}_\cK$ and a map from $N$ into $M$ with image $X_1$ covering $Z$
    such that\dots for every subset $Z'\subseteq M$ some $\prec_\cK$-extension of $N$ in level $2$
    of the tree embeds into $M$, extending the original map, and covering also $Z'$.

\Proof Let first $\gamma=0$. Then we have either $n=0$ in which case trivially $M\models
\varphi_{0,0}(=\top)$ or $n>0$. In the latter case $\varphi_{N,0,n}:=\bigwedge {\rm
Diag}^n_\kappa(N)$; if $f:N\to M$ is a $\prec_\cK$-embedding, $M$ satisfies this sentence as it
satisfies each of the formulas $\varphi(y_0,\dots y_{k-1})$ satisfied in $N$ by the images of the
$\prec_\cK$-map $f$.\\
The case $\gamma$ limit ordinal is an immediate consequence of the induction hypothesis.\\
Let now $\gamma=\beta+1$ and assume that for every $M\in \cK$, $n<\omega$, $N\in {\mathcal
S}_n$,
if $f:N\to M$ is a $\prec_\cK$-embedding then $M\models \varphi_{N,\beta,n}[\langle
f(a^*_\alpha\mid \alpha<\kappa\cdot n)\rangle]$.

Now, fix $M\in \cK$, $n<\omega$, $N\in
{\mathcal S}_n$ and $f:N\to M$ a $\cK$-embedding. We want to check that $M\models
\varphi_{N,\gamma,n}[\langle f(a^*_\alpha)\mid \alpha <\kappa\cdot n\rangle]$, i.e. we need to
verify that 
\[ M\models \forall \bar{z}_{[\kappa]}\bigvee_{\stackrel{N'\succ_\cK N}{N'\in {\mathcal
    S}_{n+1}}}
\exists \bar{x}_{=n}\Big[ \varphi_{N',\beta,n+1}(\bar{x}_n\peg\bar{x}_{=n})\wedge 
    \bigwedge_{\alpha<\kappa}\bigvee_{\delta<\kappa\cdot(n+1)}z_\alpha=x_\delta \Big]  \]
    when $\bar{x}_n$ is replaced in $M$ by $\langle f(a^*_\alpha)\mid \alpha
    <\kappa\cdot n\rangle$.

    So let $\bar{c}_{[\kappa]}\in M$. By the LST axiom, there is some $M'\prec_\cK M$ containing
    both $\bar{c}_{[\kappa]}$ and $\langle f(a^*_\alpha)\mid \alpha
    <\kappa\cdot n\rangle$, with $|M'|=\kappa$. By the isomorphism axioms there is $N'\succ_\cK
    N$, $N'\in {\mathcal S}_{n+1}$, isomorphic to $M'$ through an isomorphism $f'$ extending $f$. We may now apply the
    induction hypothesis to $N'$, $f'$: since $f':N'\to M$ is a $\prec_\cK$-embedding, we
    have that $M\models \varphi_{N',\beta,n+1}[\langle a^*_\alpha\mid \alpha
	<\kappa\cdot
    (n+1)\rangle ]$. But this enables us to conclude: $N'$ is a witness for the disjunction on
    models $\prec_\cK$-extending $N$, and
    the existential $\exists \bar{x}_{=n}$ is witnessed by $\langle a^*_\alpha\mid \alpha \in
    [\kappa\cdot n,\kappa \cdot (n+1))\rangle$. As the original $M'$ had been chosen to
    include the sequence $\bar{c}_{[\kappa]}$, the last part of the formula holds.
    \qed${}_{\mbox{Claim }\ref{inductSat}}$

In particular, when $n=0$, the previous Claim shows that if $M\in \mathcal{K}$, then
$M\models \varphi_{\gamma,0}$, for all $\gamma<\lambda^+$.

\bigskip
Now we come to the main point:

\begin{claim}\label{Mmodelsphilambda2isintheclass}
    If $M$ is a $\tau$-model and $M\models \varphi_{\lambda+1,0}$ then $M\in \cK$.
\end{claim}

\Proof The plan of this proof is as follows: we build $\mathcal G$ a set of substructures of $M$ of
cardinality $\kappa$, each of them isomorphic to a model in ${\mathcal S}_1$ and such that $M\models
\varphi_{N,\lambda,1}(\dots)$ of the elements of the substructure; we prove that $\mathcal G$ is
cofinal in $M$ (using the fact that $M\models \varphi_{\lambda+1,0}$) and a directed set. We also
prove that for elements of $\mathcal G$ being a submodel implies being a $\prec_\cK$-submodel (this
is the longest part of the proof, and requires a delicate combinatorial argument). We conclude that
$M\in \cK$, as it then ends up being the direct limit of the $\prec_\cK$-directed system $\mathcal
G$.

Let ${\mathcal G}:=\left\{ N^*\subseteq M\mid N^*\right.$ has cardinality $\kappa$ and for
some $N\in {\mathcal S}_1$ there is a bijective $f:N\to N^*$ such that $\left. M\models
\varphi_{N,\lambda,1}[\langle f(a^*_\alpha)\mid \alpha
<\kappa\rangle]\right\}$. In particular, such $f$'s are isomorphisms from $N$ to $N^*$.

We prove first
\begin{equation}
    N^*_1\subseteq N^*_2 \quad (N^*_\ell\in {\mathcal G}) \quad \mbox{ then } N^*_1\prec_\cK
    N^*_2.
    \label{modcompl}
\end{equation}

Fix $N^*_1\subseteq N^*_2$, both in $\mathcal G$.
Choose $(N^\ell_\eta,f^\ell_\eta)$ for $\ell=1,2$ and $\eta\in{\rm
ds}(\lambda):=\left\{ \nu\mid \nu\right.$ a decreasing sequence of ordinals $\left. <\lambda
\right\}$ by induction on $\ell g(\eta)$ such that
\begin{enumerate}
    \item $N^\ell_\eta\in{\mathcal S}_{\ell g(\eta)+1}$
    \item $f^\ell_\eta$ embeds $N^\ell_\eta$ into $M$: $f^\ell_\eta(N^\ell_\eta)\subseteq M$
    \item $M\models \varphi_{N^\ell_\eta,{\rm last}(\eta),\ell g(\eta)+1}[\langle
	f^\ell_\eta(a^*_\alpha\mid \alpha <\kappa\cdot (\ell g(\eta)+1))\rangle]$ where ${\rm
	last}(\langle \rangle)=\lambda$, ${\rm last}(\nu\peg\langle \alpha\rangle)=\alpha$
    \item if $\nu\triangleleft \eta$ then $N^\ell_\nu\prec_\cK N^\ell_\eta$ and $f^\ell_\nu\subseteq
	f^\ell_\eta$
    \item $f^\ell_{\langle \rangle}(N^\ell_{\langle \rangle})= N^*_\ell$
    \item $f^1_\eta(N^1_\eta)\subseteq f^2_\eta(N^2_\eta)$ and $\nu\triangleleft \eta\Rightarrow
	f^2_\nu(N^2_\nu)\subseteq
	f^1_\eta(N^1_\eta)$.
\end{enumerate}

The induction: if $\ell g(\eta)=0$ let $f^\ell_\eta=f^\ell_{\langle \rangle}$ be a one-to-one
function from $\langle a^*_\alpha\mid \alpha <\kappa\rangle $ onto $N^*_\ell$; as $\| N^*_\ell\|=\kappa$
there is a model $N^\ell_\eta$ with universe $\langle a^*_\alpha\mid \alpha
<\kappa\rangle$ such that
$f^\ell_\eta$ is an isomorphism from $N^\ell_\eta$ onto $N^*_\ell$. Since ${\rm
last}(\langle \rangle)=\lambda$ and by definition of ${\mathcal G}$ we have $M\models
\varphi_{N^\ell_{\langle \rangle},\lambda,1}[f^\ell_{\langle \rangle}(a^*_\alpha)\mid \alpha
<\kappa] $, this choice satisfies the relevant clauses (1, 2, 3, 5 and the first part of 6).

If $\ell g(\eta)=n=m+1$ we first choose $(f^1_\eta,N^1_\eta)$.
From the inductive definition of $\varphi_{N^1_{\eta\restriction m,{\rm
last}(\eta\restriction m),m}}$ with $\bar{z}_{[\kappa]}$ an enumeration of $\langle
f^2_{\eta\restriction m}(a^*_\alpha)\mid \alpha < \kappa\cdot m\rangle$, the sequence $\bar{x}_{=n}$
gives us the map $f^1_\eta$, with domain $N^1_\eta$ (a witness of the disjunction in the formula),
and $N^1_\eta\supseteq N^2_{\eta\restriction m}$. (While doing this, we make sure the new function
$f^1_\eta\supseteq f^1_{\eta\restriction m}$.)

Now to choose $(f^2_\eta,N^2_\eta)$ we use a symmetric argument and the inductive definition of
$\varphi_{N^2_{\eta\restriction m,{\rm last}(\eta\restriction m),m}}$ with
$\bar{z}_{[\kappa]}$ enumerating $\langle
f^1_{\eta\restriction n}(a^*_\alpha)\mid \alpha < \kappa\cdot n\rangle$; as before, the
sequence $\bar{x}_{=n}$ gives us the map $f^2_\eta$, with domain $N^2_\eta$. Again we make sure
$f^2_\eta\supseteq f^2_{\eta\restriction m}$.

In both construction steps the model obtained is a $\prec_\cK$-extension, since it is given by the
disjunction inside the formula $\varphi_{N^\ell_{\eta\restriction m,{\rm
last}(\eta\restriction m),m}}$.

This finishes the inductive construction of the well-founded tree of models and functions
$(N^\ell_\eta,f^\ell_\eta)_{\eta\in {\rm ds}(\lambda)}$.

Let us now check why having carried the induction suffices.

We apply a partition theorem on well founded trees due to Komjáth and Shelah~\cite{KoSh:796}.
In~\cite{GhSh:909}, Gruenhut and Shelah provide the following
useful form.

\begin{theorem}[Komjáth-Shelah, \cite{KoSh:796}]
    Let $\alpha$ be an ordinal and $\mu$ a cardinal. Set $\nu=\left(
    |\alpha|^{\mu^{\aleph_0}} \right)^+$ and let $F({\rm ds}(\nu^+))\to \mu$ be a colouring of
    the tree of strictly decreasing sequences of ordinals $<\lambda$. Then there is an embedding
    $\varphi:{\rm ds}(\alpha)\to {\rm ds}(\nu^+)$ and a function $c:\omega\to\mu$ such
    that for every $\eta\in{\rm ds}(\alpha)$ of length $n+1$
    \[ F(\varphi(\eta))=c(n).\]
    \label{kosh796}
\end{theorem}


In our case, the number of colors $\mu$ will be $\kappa^{|\tau|+\kappa}=2^\kappa$. So, the corresponding
$\nu$ is $\left(|\alpha|^{\mu^{\aleph_0}} \right)^+=\left(|\alpha|^{(2^\kappa)^{\aleph_0}} \right)^+
= \left( |\alpha|^{2^\kappa} \right)^+ = \beth_2(\kappa)^+$ hence $\nu^+=\beth_2(\kappa)^{+
+}=\lambda$. Our coloring (given 
by the choice of the models $N^\ell_\eta$ and equivalence classes\footnote{We call a
    composition of the form $g_{\eta}:=(f^2_{\eta})^{-1}\circ f^1_{\eta}:N^1_\eta \to
    N^2_{\eta}$ a \textbf{color}: there are at most as many colors as maps from models of the
    form $N^1_\eta$ 
    into the corresponding models $N^2_\eta$, but these models are all of cardinality $\kappa$
    and fixed in advance. The number of colors is therefore bounded by
    $\mu=\kappa^{|\tau|+\kappa}=2^{\kappa}$.
     The authors thank Will Boney
for requesting a clarification of this important point.} of maps $f^\ell_\eta$ for $\eta\in {\rm ds}(\lambda)$) is
therefore a mapping
\[ F:{\rm ds}(\lambda)\to \mu\]
and the partition theorem provides a sequence
$(\eta_n)_{n<\omega}$,
$\eta_n\in {\rm ds}(\alpha)$ such that: 

\[ k\leq m\leq n,\ell\in\left\{ 1,2 \right\} \Rightarrow N^\ell_{\eta_m\restriction
k}=N^\ell_{\eta_n\restriction k}.\]

We therefore obtain $(N^\ell_k,g^\ell_{k,n})_{k\leq n}$ such that

\begin{itemize}
    \item $N^1_k\subseteq N^2_k\subseteq N^1_{k+1}$ and
    \item $g^\ell_{k,n}$ is an isomorphism from $N^\ell_k$ onto
	$N^\ell_{\eta_n\restriction k}$.
\end{itemize}

Hence $N^\ell_n\prec_\cK N^\ell_{n+1}$ and so $\langle N^\ell_n\mid
n<\omega\rangle$ is $\prec_\cK$-increasing. Let $N_\ell:=\bigcup_nN^\ell_n$. Then clearly
$N_1=N_2$; call this model $N$. Since we then have $N^1_n\prec_\cK N$, $N^2_n\prec_\cK N$ and
$N^1_n\subseteq N^2_n$ by the coherence axiom for a.e.c.'s we have that $N^1_n\prec_\cK N^2_n$. In
particular, when $n=0$ we get that $N^*_1\prec_\cK N^*_2$.

Finally, we also have that
\begin{equation}
    {\mathcal G}\mbox{ is cofinal in }[M]^{\leq \kappa},
    \label{NcofinalinMkappa}
\end{equation}
as  $M\models \varphi_{\lambda+1,0}$ and the definition of the sentence $\varphi_{\lambda,0}$ says
that every $Z\subseteq M$ can be covered by some $N^*$ of cardinality $\kappa$ isomorphic to
\emph{some} $N\in {\mathcal S}_1$ such that $M\models \varphi_{N,\lambda,1}(\langle
f(a^*_\alpha)\mid \alpha <\kappa\rangle )$\dots but this means $N^*\in {\mathcal G}$.
Also, $\mathcal G$ is a directed system.

Finally, putting together~(\ref{modcompl}) and~(\ref{NcofinalinMkappa}), we conclude that every
$\tau$-model $M$ such that $M\models \varphi_{\lambda+1,0}$ must be in the class: $M=\bigcup
{\mathcal G}$, and ${\mathcal G}$ is a $\prec_\cK$-directed system. Since $\cK$ is an a.e.c, the
limit of this $\prec_\cK$-directed system must be an element of $\cK$, therefore $M\in \cK$.
\qed${}_{{\rm Lemma}~\ref{Mmodelsphilambda2isintheclass}}$

Lastly, we complete the proof of
Theorem~\ref{aecDefinableInfinitary}: Claims~\ref{inductSat}
and~\ref{Mmodelsphilambda2isintheclass} provide the definability in the class, as
clearly $\varphi_{\gamma,0}\in {\mathbb L}_{\lambda^+,\kappa^+}(\tau_\cK)$.

\qed${}_{ {\rm Theorem}~\ref{aecDefinableInfinitary}}$

%
%
%
%
%

\section{Strong embeddings and definability}


We now focus on the relation $\prec_\cK$ of our a.e.c. $\cK$: we characterize it in ${\mathbb
L}_{\lambda^+,\kappa^+}(\tau)$. We prove a syntactic criterion for being a $\prec_\cK$-substructure
(given that we already have that $M_1\subseteq M_2$) \emph{in terms of satisfiability in} $M_2$ of
certain formulas on tuples from $M_1$. This may be regarded as a very strong analog of a
``Tarski-Vaught'' criterion for a.e.c.'s.

It is worth mentioning we will continue using in a crucial way both the canonical tree ${\mathcal
S}_\cK$ of our a.e.c., and the partition theorem on well-founded trees.

\begin{theorem}
    Let $\cK$ be an a.e.c., $\tau=\tau(\cK)\leq \kappa=LST(\cK)$, $\lambda=\beth_2(\kappa)^{++}$.
    \underline{Then,} given $\tau$-models $M_1\subseteq M_2$, the following are equivalent:
    \begin{enumerate}[label=(\Alph*)]
    \item $M_1\prec_\cK M_2$
    \item if $\bar{a}_\ell \in {}^{\kappa\geq}(M_\ell)$ for $\ell=1,2$ and $\gamma <\lambda$ then
	there are $\bar{b}_\ell$, $N_\ell$ and $f_\ell$ for $\ell=1,2$ such that:\\
	for $\ell=1,2$,
	\begin{enumerate}
	    \item $\bar{b}_\ell\in {}^{\kappa\geq}(M_\ell)$ and $N_\ell\in {\mathcal S}_\ell$
	    \item ${\rm Rang}(\bar{a}_\ell)\subseteq {\rm Rang}(\bar{b}_\ell)$
	    \item $f_\ell$ is an isomorphism from $N_\ell$ onto $M_\ell\restriction {\rm
		Rang}(\bar{b}_\ell)$
	    \item ${\rm Rang}(\bar{b}_1)\subseteq {\rm Rang}(\bar{b}_2)$
	    \item $N_1\subseteq N_2$
	    \item $M_\ell\models \varphi_{N_\ell,\gamma,\ell}[\langle f_\ell(a^*_\alpha)\mid \alpha
		<\kappa\cdot \ell\rangle]$.
	\end{enumerate}
    \item if $\bar{a}\in {}^{\kappa\geq}(M_1)$ then there are $\bar{b}$, $N$ and $f$ such that
	\begin{enumerate}
	    \item $\bar{b}\in {}^{\kappa\geq}(M_1)$ and $N\in {\mathcal S}_1$
	    \item ${\rm Rang}(\bar{a})\subseteq {\rm Rang}(\bar{b})$
	    \item $f$ is an isomorphism from $N$ onto $M_1\restriction {\rm Rang}(\bar{b})$
	    \item $M_2\models \varphi_{N,\lambda+1,1}[\langle f(a^*_\alpha)\mid\alpha<\kappa)\rangle]$.
	\end{enumerate}
    \end{enumerate}
    \label{newTVlike}
\end{theorem}

\Proof
\underline{\emph{(A)}$\Rightarrow$ \emph{(B)}}: Let $\bar{a}_\ell \in {}^{\kappa\geq}(M_\ell)$ for
$\ell=1,2$ and let $\gamma <\lambda$. Choose first $N_1^*\prec_\cK M_1$ of cardinality $\leq \kappa$
including ${\rm Rang}(\bar{a}_1)$ and next, choose $N_2^*\prec_\cK M_2$ including $N^*_1 \cup
\bar{a}_2$, of cardinality $\kappa$. Let $\bar{b_\ell}$ enumerate
$N^*_\ell$ and let $(N_1,f_1,N_2,f_2)$ be such that
\begin{enumerate}
    \item $N_1\in{\mathcal S}_1$, $N_2\in{\mathcal S}_2$, $N_1\subseteq N_2$ and
    \item $f_\ell$ is an isomorphism from $N_\ell$ onto $N_\ell^*$ for $\ell=1,2$.
\end{enumerate}

This is possible: since $M_1\prec_\cK M_2$ and $N^*_\ell\prec_\cK M_\ell$ for $\ell=1,2$, we also
have that $N^*_1\prec_\cK N^*_2$. Therefore there are corresponding models $N_1\subseteq N_2$ in the
canonical tree, at levels $1$ and $2$ (as these must satisfy $N_1\prec_\cK N_2$).

We then have that $f_\ell:N_\ell\to M_\ell$ is a $\cK$-embedding from elements $N_1$ and $N_2$ in the
canonical tree $\mathcal S$. By Claim~\ref{inductSat}, we may conclude that
\[M_1\models
\varphi_{N_1,\gamma,1}[\langle f(a^*_\alpha)\mid \alpha<\kappa\rangle]\]
and
\[M_2\models
\varphi_{N_2,\gamma,2}[\langle f(a^*_\alpha)\mid \alpha<\kappa\cdot 2\rangle],\]
for each $\gamma<\lambda$.

\noindent
\underline{\emph{(B)}$\Rightarrow$ \emph{(C)}}: let $\bar{a}\in {}^{\kappa\geq}(M_1)$. We need
$\bar{b}$, $N\in {\mathcal S}_1$ and $f:N\to M_1\restriction {\rm Rang}(\bar{b})$ such that
\begin{equation}
M_2\models \varphi_{N,\lambda+1,1}[\langle f(a^*_\alpha\mid \alpha<\kappa)\rangle].
    \label{fullSatisf}
\end{equation}
(B) provides a model $N=N_1\in {\mathcal S}_1$ and elements $\bar{b}=\bar{b}_1$, as well as an
isomorphism $f:N\to {\rm Rang}(\bar{b})$. We now check that (B) also implies~\ref{fullSatisf}.

Recall the definition of $\varphi_{N,\lambda+1,1}$ (as applied to $[\langle f(a^*_\alpha\mid
\alpha<\kappa)\rangle]$).
This formula holds in $M_2$ if for every $\bar{c}_{[\kappa]}$ (of size $\kappa$) in $M_2$,
\textbf{for some} $\prec_\cK$-extension $N'$ of $N$ in ${\mathcal S}_2$ we have that
\begin{equation}
    M_2\models \exists \bar{x}_{=2}\varphi_{N',\lambda,2}[\langle f(a^*_\alpha\mid
    \alpha<\kappa)\rangle\peg \bar{x}_{=2}]
    \label{extensionSatisf}
\end{equation}
and the elements $\bar{c}_{[\kappa]}$ are ``covered'' by the list of elements (of length
$\kappa\cdot 2$) $\langle f(a^*_\alpha\mid
\alpha<\kappa)\rangle\peg \bar{x}_{=2}$. But the remaining part of clause (B) provides just this: there is
\emph{some} $N'=N_2\in {\mathcal S}_2$, extending $N=N_1$ such that for each $\gamma<\lambda$, and
an isomorphism $f'$ from $N'$ into some $\prec_\cK$-submodel $N^*$ of $M_2$ containing ${\rm
Rang}(\bar{c}_{[\kappa]})$ such that
$M_2\models \varphi_{N',\gamma,2}[\langle f'(a^*_\alpha\mid \alpha<\kappa\cdot 2)\rangle]$. The
submodel $N'$ witnesses the disjunction on models and $\langle f'(a^*_\alpha)\mid
\alpha\in[\kappa,\kappa\cdot 2)\rangle$ witnesses the existential $\bar{x}_{=2}$. 

\noindent
\underline{(\emph{C})$\Rightarrow$ (\emph{A})}: assuming (C) means that for every
$\kappa$-tuple $\bar{a}$ from $M_1$ there are a model $N\in {\mathcal S}_1$, a
$\kappa$-tuple $\bar{b}$ from $M_1$ containing $\bar{a}$ and an isomorphism from $N$ onto
$M_1\restriction {\rm Rang}(\bar{b})$ such that
\[ M_2\models \varphi_{N,\lambda+1,1}[\langle f(a^*_\alpha)\mid \alpha<\kappa\rangle ].\]
This means that for each $\bar{c}$ included in $M_2$ (of length $\kappa$) there are some extension
$N'$ of $N$ with $N'\in {\mathcal S}_2$ and some $\bar{d}$ included in $M_2$, of length
$\kappa$,
such that
\[ M_2\models \varphi_{N',\lambda,2}[\langle f(a^*_\alpha)\mid
\alpha<\kappa\rangle\peg \bar{d}]\]
    and such that ${\rm Rang}(\bar{c})\subseteq {\rm Rang}([\langle f(a^*_\alpha)\rangle]\peg
\bar{d}])$.

Consider first the family
\[ {\mathcal G}_1:=\left\{ N^*_1\subseteq M_1\mid \exists N_1\in {\mathcal S}_1 \exists
    f:N_1\stackrel{\approx}{\rightarrow}N^*_1 \big[
M_2\models\varphi_{N_1,\lambda+1,1}(f(a^*_\alpha)_{\alpha<\kappa}) \big] \right\};\]
by part $(d)$ of the hypothesis ${\mathcal G}_1$ is a directed family, cofinal in $M_1$.

Now fix $N^*_1\in {\mathcal G}_1$ and let
\[ {\mathcal G}_{N^*_1,2}:=\Big\{ N^*_2\subseteq M_2\mid N^*_1\prec_\cK N^*_2\mbox{ and} \qquad \qquad
    \qquad\qquad\]
    \[ \qquad \qquad \exists N_2\in {\mathcal S}_2 \exists
    f_2:N_2\stackrel{\approx}{\rightarrow}N^*_2 \big[
M_2\models\varphi_{N_2,\lambda,2}(f_2(a^*_\alpha)_{\alpha<\kappa\cdot 2}) \big] \Big\}.\]
Now build a tree of models as in the proof of Claim~\ref{Mmodelsphilambda2isintheclass} inside
${\mathcal S}_\cK$, indexed by ${\rm ds}(\lambda)$, and use the partition theorem on well-founded
trees to conclude that
\begin{equation}
N^*_{2,1}\subseteq N^*_{2,2}, N^*_{2,\ell}\in {\mathcal G}_{N^*_1,2} (\ell=1,2) \Rightarrow
N^*_{2,1}\prec_\cK N^*_{2,2} .
    \label{embedisstrong}
\end{equation}
Now, one of the consequences of $M_2\models \varphi_{N_1,\lambda+1,1}(f(a^*_\alpha)_{\alpha
<\kappa})$ (for the model $N_1$ in ${\mathcal S}_1$ corresponding to $N^*_1$ and for the map $f$) is
precisely that ${\mathcal G}_{N^*_1,2}$ is cofinal in $M_2$ and a directed family,
$\prec_\cK$-directed also, by~(\ref{embedisstrong}). Therefore, by the union axiom of a.e.c.'s we may
conclude that $N^*_2\prec_\cK M_2$; since we also had $N^*_1\prec_\cK N^*_2$, we have that
$N^*_1\prec_\cK M_2$.

Since $N^*_1$ was an arbitrary member of ${\mathcal G}_1$, we may conclude that all members of
${\mathcal G}_1$ are $\prec_\cK$-elementary in $M_2$. By another application of the partition
relation, the family ${\mathcal G}_1$ also has the property that $M^*_{1,1}\subseteq M^*_{1,2}$ in
the family implies $M^*_{1,1}\prec_\cK M^*_{1,2}$. So, applying again the union axiom, we may
conclude that $M_1=\bigcup {\mathcal G}_1\prec_\cK M_2$.

%
%
%
%
%

    \qed${}_{{\rm Theorem}~\ref{newTVlike}}$

    The previous criterion for $M_1\prec_\cK M_2$, given $M_1\subseteq M_2$, is admittedly quite
    sophisticated compared with the classical Tarski-Vaught criterion for elementarity in first
    order logic. There are, however, some interesting parallels.

    \begin{itemize}
	\item In part (C) of our criterion, we only evaluate the formula at the ``large model''
	    $M_2$. This is one of the crucial aspects of the Tarski-Vaught criterion, as it allows
	    construction ``from below'' of elementary submodels.
	\item The aspect of our criterion that is definitely less within reach is a version of
	    ``capturing existential formulas.'' We are in a sense exactly doing \emph{that} but in
	    the more complex world of a.e.c.'s. Satisfying a formula of the form
	    $\varphi_{N,\lambda+1,1}(\dots)$ at a subset of elements of the \emph{small} model
	    $M_1$, when parsing the formula, in a way reflects the possibility of being able to
	    realize, according to $M_2$, all ``possible extensions'' of small models, reflecting
	    them correctly to $M_1$. The partition relation on well-founded trees of course ends up
	    being the key in our case.
    \end{itemize}

    \section{Around the logic of an a.e.c.}

The logic usually called $L^1_\kappa$ from Shelah's paper~\cite{Sh:797} satisfies Interpolation and a weak form of compactness: strong undefinability of well-order. Furthermore, it satisfies
a Lindström-like maximality theorem for these properties (as well as union of $\omega$-chains of
models). The logic $L^1_\kappa$, however, has a non-algorithmic syntax (sentences are unions of
equivalence classes of structures under a relation defined based on a ``delayed Ehrenfeucht-Fraïssé
game''). Väänänen and
Villaveces~\cite{VaVi} have produced a logic with a clearly defined (and relatively symple) syntax,
whose $\Delta$-closure (a notion appearing first in~\cite{MShS:47}) is $L^1_\kappa$, and which
satisfies several of the good properties of that logic (of course, strong undefinability of
well-order but also closure under unions of chains). Also, D\v{z}amonja and Väänänen
have proved in~\cite{DzVa} that \emph{chain logic} shares many properties with Shelah's logic
$L^1_\kappa$, with an important semantic difference (the notion of a model); they have also provided
careful comparisons between the two logics.

All of these logics are close to our constructions in this paper: the sentence
$\varphi_{\lambda+1,0}$ belongs to ${\mathbb L}_{\lambda^+,\kappa^+}$ and $L^1_\mu$ lies in
between two logics of the form $L_{\mu,\aleph_0}$ and $L_{\mu,\mu}$. Our sentence
$\varphi_{\lambda+1,0}$ belongs to $L^1_\mu$. However, it is not clear if this is the minimal logic
for which this is the case.

The question of which is the minimal logic capturing an a.e.c. remains still partially
open. Our theorems in this article provide a substantial advance in this direction.

\end{document}